# The Askey-Wilson Polynomials and $q$-Sturm-Liouville Problems[*]

B. Malcolm Brown, W. Desmond Evans and Mourad E. H. Ismail

August 14, 1994

## Abstract

We find the adjoint of the Askey-Wilson divided difference operator with respect to the inner product on $L^2(-1, 1, (1-x^2)^{-1/2}dx)$ defined as a Cauchy principal value and show that the Askey-Wilson polynomials are solutions of a $q$-Sturm-Liouville problem. From these facts we deduce various properties of the polynomials in a simple and straightforward way. We also provide an operator theoretic description of the Askey-Wilson operator.

**Running title**: $q$-Sturm-Liouville Problems.

**1. Introduction**. Askey and Wilson [5] introduced a remarkable system of orthogonal polynomials $\{p_n(x; a, b, c, d|q)\}$ depending on four parameters $a, b, c, d$ and a base parameter $q$ which became known as the Askey-Wilson polynomials. They generalize the classical polynomials of Jacobi, Hermite and Laguerre, [3]. Askey and Wilson [5] also introduced a divided difference operator $\mathcal{D}_q$, whose action on the Askey-Wilson polynomials mimics the action of $\frac{d}{dx}$ on Jacobi polynomials.

Our first step in this paper is to compute $\mathcal{D}_q^*$, the adjoint of $\mathcal{D}_q$ with respect to the $L^2(-1, 1, (1-x^2)^{-1/2}dx)$ inner product defined as a Cauchy principal value. This result and the raising operator for the Askey-Wilson polynomials are then used to investigate the properties of the polynomials as solutions of a $q$-Sturm-Liouville equation. It is shown how a number of properties of the Askey-Wilson polynomials, like their orthogonality and a Rodrigues-type formula, can be deduced in a way which is much simpler than those found in the literature.

---

[*]Research partially supported by NSF grant DMS 9203659, SERC grant GR/J54604 and NATO grant 1382/93/JARC-50.

In this work we mostly follow the terminology of [3] and [7]. We will always assume

$$0 < q < 1.$$

We first remind the reader of the notations to be used. A $q$-shifted factorial is defined by [7]

$$(1.1) \quad (a;q)_0 := 1, \quad (a;q)_n := \prod_{k=1}^{n}(1 - aq^{k-1}), \quad n = 1, 2, \cdots, \infty,$$

and more generally [7]

$$(1.2) \quad (a_1, \ldots, a_k; q)_n = \prod_{j=1}^{k} \prod_{r=1}^{n}(1 - a_j q^{r-1}), \quad n = 0, 1, 2, \ldots, \infty.$$

A basic hypergeometric series is

$$(1.3) \quad {}_{r+1}\phi_r \left( \begin{array}{c} a_1, \ldots, a_{r+1} \\ b_1, \ldots, b_r \end{array} \bigg| q, z \right) = \sum_{n=0}^{\infty} \frac{(a_1, \ldots, a_{r+1}; q)_n}{(q, b_1, \ldots, b_r; q)_n} z^n.$$

The Askey-Wilson polynomials are [5], [7, (7.5.2)]

$$(1.4) \quad p_n(\cos\theta; a, b, c, d|q) = (ab, ac, ad; q)_n a^{-n} {}_4\phi_3 \left( \begin{array}{c} q^{-n}, abcdq^{n-1}, ae^{i\theta}, ae^{-i\theta} \\ ab, \quad ac, \quad ad \end{array} \bigg| q, z \right),$$

for $(n = 0, 1, 2, \ldots)$. When $a, b, c, d \in (-1, 1)$ the Askey-Wilson polynomials satisfy the orthogonality relation

$$(1.5) \quad \int_{-1}^{1} w(x; a, b, c, d) \, p_n(x; a, b, c, d|q) \, p_m(x; a, b, c, d|q) \, dx = \zeta_n(a, b, c, d) \, \delta_{m,n},$$

where

$$(1.6) \quad w(x; a, b, c, d) := \frac{(1 - x^2)^{-1/2}(e^{2i\theta}, e^{-2i\theta}; q)_\infty}{(ae^{i\theta}, ae^{-i\theta}, be^{i\theta}, be^{-i\theta}, ce^{i\theta}, ce^{-i\theta}, de^{i\theta}, de^{-i\theta}; q)_\infty}, \quad x := \cos\theta,$$

and

$$(1.7) \quad \zeta_n(a, b, c, d) = \frac{2\pi (abcdq^{2n}; q)_\infty (abcdq^{n-1}; q)_n}{(q^{n+1}, abq^n, acq^n, adq^n, bcq^n, bdq^n, cdq^n; q)_\infty}.$$

Note that

$$w(x; 1, -1, q^{1/2}, -q^{1/2}) = (1 - x^2)^{-1/2}, \quad w(x; q^{1/2}, -q^{1/2}, q, -q) = 4(1 - x^2)^{1/2}.$$

Therefore the special Askey-Wilson polynomials $\{p_n(x; 1, -1, q^{1/2}, -q^{1/2}|q)\}$ and $\{p_n(x; q^{1/2}, -q^{1/2}, q, -q|q)\}$ are multiples of the Chebyshev polynomials of the first and second kinds, respectively. The Jacobi polynomials and hence the Hermite and Laguerre polynomials are limits of special Askey-Wilson polynomials. Indeed

$$P_n^{\alpha,\beta}(x) = \lim_{q \to 1} \frac{p_n(x; q^{(2\alpha+1)/4}, q^{(2\alpha+3)/4}, -q^{(2\beta+1)/4}, q^{(2\beta+3)/4}|q)}{(q, -q^{(\alpha+\beta+1)/2}, -q^{(\alpha+\beta+2)/2}; q)_n}.$$

The Askey-Wilson memoir [5] contains a list of special and limiting cases of the Askey-Wilson polynomials.

The inner product which we shall use initially in this work is the inner product associated with the Chebyshev weight $(1 - x^2)^{-1/2}$ on $(-1, 1)$, namely

$$(1.8) \quad < f, g > := \int_{-1}^{1} f(x) \, \overline{g(x)} \, \frac{dx}{\sqrt{1 - x^2}}.$$

Given a function $f$ defined on $(-1, 1)$ we set $\breve{f}(e^{i\theta}) := f(x), \; x = \cos\theta$, that is

$$(1.9) \quad \breve{f}(z) = f((z + 1/z)/2), \, z = e^{i\theta}.$$

In this notation the Askey-Wilson finite difference operator $\mathcal{D}_q$ is defined by

$$(1.10) \quad (\mathcal{D}_q f)(x) := \frac{\breve{f}(q^{1/2}e^{i\theta}) - \breve{f}(q^{-1/2}e^{i\theta})}{(q^{1/2} - q^{-1/2}) \, i \, \sin\theta}, \quad x = \cos\theta.$$

This requires $\breve{f}(z)$ to be defined for $|q^{\pm 1/2}z| = 1$ as well as for $|z| = 1$. In particular $\mathcal{D}_q$ is well defined on $H_{\frac{1}{2}}$, where

$$(1.11) \quad H_\nu := \{f : f((z + 1/z)/2) \text{ is analytic for } q^\nu \leq |z| \leq q^{-\nu}\}.$$

For an interesting way in which the Askey-Wilson operator arises naturally from divided difference operators we refer the reader to [15]. Askey [4] used the Askey-Wilson operator to give a new proof of the connection coefficient formula for the continuous $q$-ultraspherical polynomials. Ismail [9] proved a Taylor-like series for polynomials, with $\frac{d}{dx}$ replaced by $\mathcal{D}_q$, and used it to give proofs of several summation theorems and transformation formulas for basic hypergeometric functions.

This work is a contribution to the study of Askey-Wilson polynomials and $q$-Sturm-Liouville problems. We shall show that the orthogonality relation (1.5) is a special case of the orthogonality of different eigenspaces of a $q$-Sturm-Liouville problem. The use of the inner product (1.8) in the context of the Askey-Wilson polynomials is also new. We shall only assume the explicit representaion (1.4).

In Section 2 we shall prove that, for $f \in H_{1/2}$ and $g \in H_{1/2}$ the identity

$$(1.12) \quad <\mathcal{D}_q f, g> \quad = \frac{\pi \sqrt{q}}{1-q} \left[ f(\frac{1}{2}(q^{1/2}+q^{-1/2}))\overline{g(1)} - f(-\frac{1}{2}(q^{1/2}+q^{-1/2}))\overline{g(-1)} \right]$$
$$- <f, \sqrt{1-x^2}\, \mathcal{D}_q(g(x)\,(1-x^2)^{-1/2})>$$

holds. Formula (1.12) is the exact analogue of integration by parts and suggests that the formal adjoint of $\mathcal{D}_q$ with respect to $<\cdot,\cdot>$ is

$$(1.13) \quad \mathcal{D}_q^* := -\sqrt{1-x^2}\, \mathcal{D}_q \, \frac{1}{\sqrt{1-x^2}}.$$

In Section 3 we shall consider the $q$-Sturm-Liouville problem

$$(1.14) \quad \frac{1}{w(x)}\mathcal{D}_q(\,p(x)\,\mathcal{D}_q\,y(x)) = \lambda\,y(x),$$

and prove that when $w(x) > 0$ almost everywhere on $(-1,1)$ the appropriate solutions of (1.14) with distinct $\lambda$'s are orthogonal in $L^2(-1,1,(1-x^2)^{-1/2}dx)$. We shall prove in Section 3 that a particular case of (1.14) is satisfied by the Askey-Wilson polynomials, namely

$$(1.15) \quad \frac{1}{w(x;a,b,c,d)}\mathcal{D}_q\left(w(x;aq^{1/2},bq^{1/2},cq^{1/2},dq^{1/2})\,\mathcal{D}_q\,p_n(x;a,b,c,d|q)\right) = \lambda_n\,p_n(x;a,b,c,d|q),$$

with

$$(1.16) \quad \lambda_n := \frac{4q}{(1-q)^2}(1-q^{-n})(1-abcdq^{n-1}), \quad n = 0,1,2,\cdots.$$

In addition, we shall prove in Theorem 3.4 that

$$(1.17) \quad \frac{1}{w(x;a,b,c,d)}\mathcal{D}_q\left[w(x;aq^{1/2},bq^{1/2},cq^{1/2},dq^{1/2})\,\mathcal{D}_q\,f(x)\right] = \lambda\,f(x),$$

has a polynomial solution if and only if $\lambda = \lambda_n$, in which case the polynomial solution is $p_n(x;a,b,c,d)$. Note that in the case of the Chebyshev polynomials of the first kind $w(x;1,-1,q^{1/2},-q^{-1/2}) = (1-x^2)^{-1/2}$ and $w(x;q^{1/2},-q^{1/2},q,-q) = 4(1-x^2)^{1/2}$ and (1.15) becomes

$$\sqrt{1-x^2}\mathcal{D}_q\left[\sqrt{1-x^2}\mathcal{D}_q T_n(x)\right] = -\frac{(q^{-n/2}-q^{n/2})^2}{(q^{-1/2}-q^{1/2})^2}T_n(x).$$

Our approach to the study of the Askey-Wilson polynomials differs from those of others, [5], [6], [11], [12], [13], [16], mainly in the use we make of (1.12) and (1.15) and in our evaluation of the $L^2$ norms of the Askey-Wilson polynomials. Kalnins and Miller [13] used a weighted space approach and do not seem to have noticed the Hilbert space connection found here.

Section 4 contains several remarks, some of which are of a pedagogical nature, which we hope will clarify certain points pertaining to our approach. Section 4 also contains a simple evaluation of the connection coefficients in

$$(1.18) \quad p_n(x; a, b, aq^{1/2}, bq^{1/2}|q) = \sum_{k=0}^{n} c_{n,j}(a,b) p_k(x; aq^{1/2}, bq^{1/2}, aq, bq|q).$$

The two-parameter family of polynomials $\{p_n(x; a, b, aq^{1/2}, bq^{1/2}|q)\}$ are called the continuous $q$-Jacobi polynomials and are $q$-analogues of the Jacobi polynomials, [5], [7]. Ismail, Rahman and Zhang [10] evaluated the connection coefficients $\{c_{n,j}(a,b)\}_{j=0}^{n}$ using a long and complicated procedure. Here the $c_{n,j}$'s are determined in a simple way from (1.12). The solution of the connection coefficient problem (1.18) was needed in [8] to find a tridiagonal representation of an inverse to $\mathcal{D}_q$ on certain weighted $L^2$ spaces.

In Section 5 we describe the problem in a Hilbert space setting and examine the Friedrichs extension of a positive operator generated by

$$(1.19) \quad (Mf)(x) := -\frac{1}{w(x)} \mathcal{D}_q \left(p \mathcal{D}_q f\right)(x).$$

**2. $q$-Integration By Parts**. In this section we shall give a proof of (1.12) which is the $q$-analogue of integration by parts for the Askey-Wilson operator $\mathcal{D}_q$.

**Proof of (1.12)**. Clearly

$$(2.1) \quad <\mathcal{D}_q f, g> = \int_0^{\pi} \frac{\check{f}(q^{1/2} e^{i\theta}) - \check{f}(q^{-1/2} e^{i\theta})}{(q^{1/2} - q^{-1/2}) i \sin\theta} \overline{g(\cos\theta)} d\theta.$$

The integral in (2.1) is defined as a Cauchy principal value and we consider

$$(2.2) \quad I_\epsilon := \int_\epsilon^{\pi-\epsilon} \frac{\check{f}(q^{1/2} e^{i\theta}) - \check{f}(q^{-1/2} e^{i\theta})}{(q^{1/2} - q^{-1/2}) i \sin\theta} \overline{g(\cos\theta)} d\theta.$$

The observation

$$\check{f}(q^{-1/2} e^{-i\theta}) = f((q^{-1/2} e^{-i\theta} + q^{1/2} e^{i\theta})/2) = \check{f}(q^{1/2} e^{i\theta})$$

leads to

$$I_\epsilon = \int_\epsilon^{\pi-\epsilon} \frac{\breve{f}(q^{1/2}e^{i\theta})}{i(q^{1/2}-q^{-1/2})} \frac{\overline{g(\cos\theta)}}{\sin\theta} d\theta - \int_{-\epsilon}^{\epsilon-\pi} \frac{\breve{f}(q^{-1/2}e^{-i\theta})}{i(q^{1/2}-q^{-1/2})} \frac{\overline{g(\cos\theta)}}{\sin\theta} d\theta$$

$$= \left(\int_{-\pi+\epsilon}^{-\epsilon} + \int_\epsilon^{\pi-\epsilon}\right) \frac{\breve{f}(q^{1/2}e^{i\theta})}{i(q^{1/2}-q^{-1/2})} \frac{\overline{g(\cos\theta)}}{\sin\theta} d\theta.$$

It is clear from Figure 1 that

$$I_\epsilon = \left(\int_C - \int_{C_\epsilon} - \int_{C'_\epsilon}\right) \frac{\breve{f}(q^{1/2}z)}{(q^{1/2}-q^{-1/2})} \frac{\overline{\breve{g}(z)}}{(z-1/z)/2} \frac{dz}{iz},$$

where $C$ is the unit circle indented at $\pm 1$ with circular arcs centered at $\pm 1$ and radii equal to $\epsilon$; and $C_\epsilon$ and $C'_\epsilon$ are the circular arcs

$$C_\epsilon = \{z : z = 1 + \epsilon e^{i\theta}, \ -(\pi-\epsilon)/2 \le \theta \le (\pi-\epsilon)/2\}$$

and

$$C'_\epsilon = \{z : z = -1 + \epsilon e^{i\theta}, \ -(\pi-\epsilon)/2 \le \theta \le (\pi-\epsilon)/2\}.$$

Both $C_\epsilon$ and $C'_\epsilon$ are positively oriented.

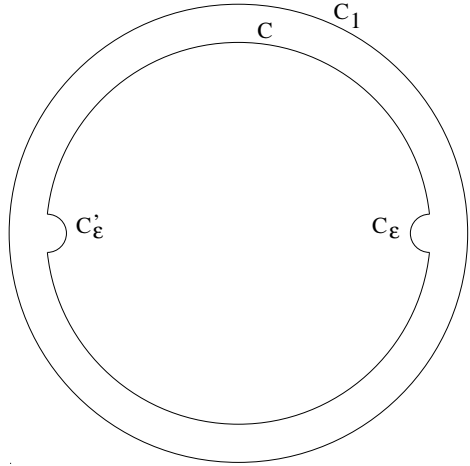

Fig 1: contour for integration

Define $\phi$ by

(2.3) $$\phi(z) := 2 \frac{\breve{f}(q^{1/2}z)}{(q^{1/2}-q^{-1/2})} \frac{\overline{\breve{g}(z)}}{iz(z-1/z)}.$$

It is clear that the residues of $\phi$ at $z = \pm 1$ are given by

$$Res\{\phi : z = 1\} = \frac{f((q^{1/2}+q^{-1/2})/2)\,\overline{g(1)}}{i(q^{1/2}-q^{-1/2})},$$

$$Res\{\phi : z = -1\} = -\frac{f(-(q^{1/2} + q^{-1/2})/2)\overline{g(-1)}}{i(q^{1/2} - q^{-1/2})}.$$

Thus

$$\lim_{\epsilon \to 0^+} \int_{C_\epsilon} \phi(z)dz = -\pi i \, Res\{\phi(z) : z = 1\}$$

and

$$\lim_{\epsilon \to 0^+} \int_{C'_\epsilon} \phi(z)dz = -\pi i \, Res\{\phi(z) : z = -1\}.$$

Let $C_1$ be the circle $|z| = q^{-1/2}$. It is clear that

$$\int_C \phi(z)dz = \int_{C_1} \phi(z)dz - 2\pi i \left[Res\{\phi(z) : z = 1\} + Res\{\phi(z) : z = -1\}\right].$$

This shows that

$$(2.4) \quad <\mathcal{D}_q f, g> \ = \ \lim_{\epsilon \to 0^+} I_\epsilon$$

$$= \frac{\sqrt{q}\pi}{1-q}[f((q^{1/2} + q^{-1/2})/2)\overline{g(1)} - f(-(q^{1/2} + q^{-1/2})/2)\overline{g(-1)}]$$

$$+ 2\int_{C_1} \frac{\check{f}(q^{1/2}z)}{(q^{1/2} - q^{-1/2})} \frac{\overline{\check{g}(z)}}{(z - 1/z)} \frac{dz}{iz}.$$

In the last integral replace $z$ by $zq^{-1/2}$ to change the integral over $C_1$ to an integral over the unit circle. The result is

$$(2.5) \quad 2\int_{C_1} \frac{\check{f}(q^{1/2}z)}{(q^{1/2} - q^{-1/2})} \frac{\overline{\check{g}(z)}}{(z - 1/z)} \frac{dz}{iz} = \int_{|z|=1} \frac{\check{f}(z)}{(q^{1/2} - q^{-1/2})} \overline{\check{h}(q^{-1/2}z)} \frac{dz}{z},$$

where

$$h(\cos\theta) = g(\cos\theta)/\sin\theta.$$

Finally in the integral over the unit circle in (2.5) we replace $z$ by $e^{i\theta}$ then write the integral as a sum of integrals over $[-\pi, 0]$ and $[0, \pi]$. In the integral over the range $[-\pi, 0]$ replace $\theta$ by $-\theta$ then combine the two integrals which are now over $[0, \pi]$. Combining this with the observation

$$\check{h}(q^{-1/2}z)|_{z=e^{-i\theta}} = -\check{h}(q^{1/2}e^{i\theta})$$

we obtain (1.12).

**3. $q$-Sturm-Liouville Problems.** In this section we study solutions to (1.14). We shall assume throughout that $p(x)/\sqrt{1-x^2} \in H_{1/2}$. Our first result is the analogue of Green's formula.

**Theorem 3.1** *For $f, g \in H_1$ we have*

(3.1) $\quad < \mathcal{D}_q\left(p(x)\mathcal{D}_q f(x)\right), \sqrt{1-x^2}g(x) > \; = \; < \sqrt{1-x^2}f(x), \mathcal{D}_q\left(p(x)\mathcal{D}_q g(x)\right) >.$

This follows from a repeated use of (1.12) with $g$ replaced by $\sqrt{1-x^2}g(x)$ and the observations that if $f \in H_{1/2}$ then so is the function $x \mapsto \sqrt{1-x^2}f(x)$ and $\mathcal{D}_q f \in H_{1/2}$ if $f \in H_1$: the boundary terms in (3.1) vanishe because of the factor $\sqrt{1-x^2}$.

We next consider eigenfunctions corresponding to different eigenvalues.

**Theorem 3.2** *Let $y_1, y_2 \in H_1$ be solutions to (1.14) with $\lambda = \lambda_1$ and $\lambda = \lambda_2$, respectively and assume $\lambda_1 \neq \lambda_2$. If $w(x) \geq 0$, $-1 \leq x \leq 1$ and $w(x) \neq 0$ almost everywhere on $[-1, 1]$, then*

(3.2) $\quad \int_{-1}^{1} w(x) y_1(x) y_2(x) dx = 0.$

*Furthermore, the eigenvalues of (1.14) are all real.*

**Proof**. It follows from (1.14) that

$$(\lambda_1 - \lambda_2)\int_{-1}^{1} w(x) y_1(x) \overline{y_2(x)} dx$$

$$= \; < \mathcal{D}_q p(x) \mathcal{D}_q y_1(x), \sqrt{1-x^2} y_2(x) > \; - \; < \sqrt{1-x^2} y_1(x), \mathcal{D}_q p(x) \mathcal{D}_q y_2(x) >$$

$$= \; 0$$

by (3.1), whence (3.2). If $\lambda_1$ is not real then $\overline{\lambda_1}$ is an eigenvalue with an eigenfunction $y_2(x) = \overline{y_1(x)}$ so the above equation holds with $\lambda_2 = \overline{\lambda_1}$. This leads to a contradiction and shows that the eigenvalues are real.

We now apply our results to give a new proof of (1.5) for the Askey-Wilson polynomials. First observe that

(3.3) $\quad \mathcal{D}_q(ae^{i\theta}, ae^{-i\theta}; q)_n = -2a\dfrac{1-q^n}{1-q}(aq^{1/2}e^{i\theta}, aq^{1/2}e^{-i\theta}; q)_{n-1}.$

Askey and Wilson used (3.3) and (1.4) to obtain

(3.4) $\quad \mathcal{D}_q p_n(x; a, b, c, d | q) = 2q^{(1-n)/2}\dfrac{(1-q^n)(1-abcdq^{n-1})}{1-q} p_{n-1}(x; aq^{1/2}, bq^{1/2}, cq^{1/2}, dq^{1/2} | q).$

Thus $\mathcal{D}_q$ is a lowering operator for the $p_n$'s since it lowers their degrees by 1. Askey and Wilson [5] used the three term recurrence relation satisfied by their polynomials, see (4.4)-(4.7), to find the raising operator

(3.5) $\quad \dfrac{2q^{(1-n)/2}}{q-1} w(x; a, b, c, d) p_n(x; a, b, c, d | q)$

$$= \mathcal{D}_q\left(w(x; aq^{1/2}, bq^{1/2}, cq^{1/2}, dq^{1/2}) p_{n-1}(x; aq^{1/2}, bq^{1/2}, cq^{1/2}, dq^{1/2}|q)\right),$$

where $w(x; a, b, c, d)$ is defined in (1.6).

We now give a direct proof of (3.5). Our proof is easier than the sketchy one of the original result in [5].

**Proof of (3.5).** It is easy to express $w(x; a, b, c, d)$ in the form

$$(3.6) \quad w(x; a, b, c, d) = \frac{2i\, e^{-i\theta}\, (qe^{2i\theta}, e^{-2i\theta}; q)_\infty}{(ae^{i\theta}, ae^{-i\theta}, be^{i\theta}, be^{-i\theta}, ce^{i\theta}, ce^{-i\theta}, de^{i\theta}, de^{-i\theta}; q)_\infty},$$

and an easy calculation using (3.6) gives

$$(3.7) \quad \frac{\mathcal{D}_q w(x; aq^{1/2}, bq^{1/2}, cq^{1/2}, dq^{1/2})}{w(x; a, b, c, d)}$$

$$= \frac{2}{q-1}[2(1 - abcd)\cos\theta - (a + b + c + d) + abc + abd + acd + bcd].$$

This leads to

$$\frac{\mathcal{D}_q\left[w(x; aq^{1/2}, bq^{1/2}, cq^{1/2}, dq^{1/2}) p_{n-1}(x; aq^{1/2}, bq^{1/2}, cq^{1/2}, dq^{1/2}|q)\right]}{w(x; a, b, c, d)}$$

$$= \frac{(abq, acq, adq; q)_{n-1}}{(a\sqrt{q})^{n-1} w(x; a, b, c, d)} \sum_{k=0}^{n-1} \frac{(q^{1-n}, abcdq^n; q)_k q^k}{(q, abq, acq, adq; q)_k}$$
$$\times \mathcal{D}_q w(x; aq^{k+1/2}, bq^{1/2}, cq^{1/2}, dq^{1/2})$$

$$= \frac{2(ab, ac, ad; q)_n}{(q-1)a^n q^{(n-1)/2}} \sum_{k=0}^{n-1} \frac{(q^{1-n}, abcdq^n; q)_k}{(q; q)_k (ab, ac, ad; q)_{k+1}} \frac{w(x; aq^k, b, c, d)}{w(x; a, b, c, d)}$$
$$\times [2aq^k \cos\theta\, (1 - abcdq^k) - aq^k(b + c + d) + a^2 q^{2k}(bc + bd + cd - 1) + bcd].$$

The term in square brackets on the right-hand side can be written as

$$(1 - abq^k)(1 - acq^k)(1 - adq^k) - (1 - aq^k e^{i\theta})(1 - aq^k e^{-i\theta})(1 - abcdq^k).$$

Therefore

$$\frac{\mathcal{D}_q\left[w(x; aq^{1/2}, bq^{1/2}, cq^{1/2}, dq^{1/2}) p_{n-1}(x; aq^{1/2}, bq^{1/2}, cq^{1/2}, dq^{1/2}|q)\right]}{w(x; a, b, c, d)}$$

$$= \frac{2(ab, ac, ad; q)_n}{(q-1)a^n q^{(n-1)/2}} \left[\sum_{k=0}^{n-1} \frac{(q^{1-n}, abcdq^n, ae^{i\theta}, ae^{-i\theta}; q)_k}{(q, ab, ac, ad; q)_k}\right.$$

$$-\sum_{k=1}^{n}\frac{(q^{1-n},abcdq^n;q)_{k-1}(1-q^k)(ae^{i\theta},ae^{-i\theta};q)_k}{(q,ab,ac,ad;q)_k}(1-abcdq^{k-1})\Bigg]$$

$$=\frac{2(ab,ac,ad;q)_n}{(q-1)a^nq^{(n-1)/2}}\sum_{k=0}^{n}\frac{(q^{1-n},abcdq^n;q)_{k-1}(ae^{i\theta},ae^{-i\theta};q)_k}{(q,ab,ac,ad;q)_k}$$

$$\times\left[(1-q^{(k-n)}(1-abcdq^{n+k-1})-(1-q^k)(1-abcdq^{k-1})\right].$$

Putting all this together establishes (3.5) since the term in square brackets is

$$q^k(1-q^{-n})(1-abcdq^{n+k-1}).$$

**Proof of (1.15).** Replace $p_{n-1}(x;aq^{1/2},bq^{1/2},cq^{1/2},dq^{1/2}|q)$ in (3.5) by its value from (3.4). Simple manipulations will establish (1.15).

It is easy to see, [5], that (3.5) implies the Rodrigues formula

$$(3.8) \quad w(x;a,b,c,d)p_n(x;a,b,c,d|q)=\left(\frac{q-1}{2}\right)^n q^{n(n-1)/4}\mathcal{D}_q^n[w(x,aq^{n/2},bq^{n/2},cq^{n/2},dq^{n/2})].$$

We now come to

**Theorem 3.3** *The Askey - Wilson polynomials satisfy (1.51).*

**Proof.** When $m\neq n$, Theorem 3.2 establishes (1.5) from (1.15), so it remains only to consider the case $m=n$. We have

$$\left(\frac{2}{q-1}\right)^n q^{n(1-n)/4}\int_{-1}^{1}p_n^2(x;a,b,c,d|q)w(x;a,b,c,d)dx$$

$$=<\mathcal{D}_q^n[w(x,aq^{n/2},bq^{n/2},cq^{n/2},dq^{n/2})],\sqrt{1-x^2}p_n(x;a,b,c,d|q)>$$

$$=(-1)^n<w(x,aq^{n/2},bq^{n/2},cq^{n/2},dq^{n/2}),\sqrt{1-x^2}\mathcal{D}_q^np_n(x;a,b,c,d|q)>,$$

where we have used (1.12) and (3.8). On applying (3.4) we obtain

$$\zeta_n(a,b,c,d)=(q,abcdq^{n-1};q)_n\,\zeta_0(aq^{n/2},bq^{n/2},cq^{n/2},dq^{n/2}).$$

But the three term recurrence relation (1.24) in [5], yields

$$(3.9) \quad \zeta_n(a,b,c,d)=\frac{(q,ab,ac,ad,bc,bd,cd;q)_n\,(abcd/q;q)_{2n}}{(abcd/q;q)_n\,(abcd;q)_{2n}}\zeta_0(a,b,c,d).$$

Thus we have established the functional equation

$$(3.10) \quad \zeta_0(a,b,c,d) = \frac{(abcd/q, abcdq^{n-1};q)_n}{(ab,ac,ad,bc,bd,cd;q)_n} \frac{(abcd;q)_{2n}}{(abcd/q;q)_{2n}} \zeta_0(aq^{n/2}, bq^{n/2}, cq^{n/2}, dq^{n/2}).$$

As $n \to \infty$ (3.10) becomes

$$(3.11) \quad \zeta_0(a,b,c,d) = \frac{(abcd;q)_\infty}{(ab,ac,ad,bc,bd,cd;q)_\infty} \zeta_0(0,0,0,0).$$

To evaluate $\zeta_0(0,0,0,0)$ we appeal to the Jacobi triple product identity, [7],

$$(3.12) \quad \sum_{-\infty}^{\infty} (-1)^k q^{k^2/2} z^k = (q, \sqrt{q}z, \sqrt{q}/z; q)_\infty,$$

and get

$$\begin{aligned}
\zeta_0(0,0,0,0) &= \int_0^\pi (e^{2i\theta}, e^{-2i\theta}; q)_\infty d\theta \\
&= \frac{1}{(q;q)_\infty} \sum_{-\infty}^{\infty} q^{k(k+1)/2} \int_0^\pi (1 - e^{2i\theta}) e^{2ki\theta} d\theta \\
&= 2\pi/(q;q)_\infty.
\end{aligned}$$

Therefore

$$(3.13) \quad \zeta_0(a,b,c,d) = \int_{-1}^{1} w(x,a,b,c,d) dx = \frac{2\pi (abcd;q)_\infty}{(q,ab,ac,ad,bc,bd,cd;q)_\infty},$$

and (3.9) shows that $\zeta_n(a,b,c,d)$ is given by (1.7). This completes the proof of Theorem 3.3.

**Theorem 3.4** *If $f$ is a polynomial solution of (1.17) of degree $n$, then $\lambda = \lambda_n$ and $f$ is a constant multiple of $p_n(x;a,b,c,d|q)$.*

**Proof.** Let

$$(3.14) \quad f(x) = \sum_{k=0}^{n} a_k (ae^{i\theta}, ae^{-i\theta}; q)_k.$$

Substitute for $f$ from (3.14) into (1.17) and use (3.3) and (3.7) to get

$$\begin{aligned}
\lambda\, w(x;a,b,c,d)\, f(x) &= \frac{2a}{q-1} \sum_{k=1}^{n} (1-q^k) a_k \mathcal{D}_q \left[ (ae^{i\theta}, ae^{-i\theta}; q)_{k-1} w(x; aq^{1/2}, bq^{1/2}, cq^{1/2}, dq^{1/2}) \right] \\
&= \frac{2a}{q-1} \sum_{k=0}^{n-1} (1-q^{k+1}) a_{k+1} \mathcal{D}_q \left[ w(x; aq^{k+1/2}, bq^{1/2}, cq^{1/2}, dq^{1/2}) \right] \\
&= \frac{4a}{(q-1)^2} \sum_{k=0}^{n-1} (1-q^{k+1}) a_{k+1} w(x; aq^k, b, c, d)
\end{aligned}$$

$$\times [2\cos\theta(1 - abcdq^k) - b - c - d + aq^k(bc + bd + cd - 1) + bcd].$$

Since
$$w(x; aq^k, b, c, d) = w(x; a, b, c, d)(ae^{i\theta}, ae^{-i\theta}; q)_k$$

we get

$$(3.15) \quad \lambda f(x) = \frac{4a}{(q-1)^2} \sum_{k=0}^{n-1} (1 - q^{k+1}) a_{k+1} (ae^{i\theta}, ae^{-i\theta}; q)_k$$

$$\times [2\cos\theta(1 - abcdq^k) - b - c - d + aq^k(bc + bd + cd - 1) + bcd].$$

As in the proof of (3.5) we write

$$(3.16) \quad aq^k[2\cos\theta(1 - abcdq^k) - b - c - d + aq^k(bc + bd + cd - 1) + bcd]$$

$$= (1 - abq^k)(1 - acq^k)(1 - adq^k) - (1 - abcdq^k)(1 - aq^k e^{i\theta})(1 - aq^k e^{-i\theta}).$$

We now use (3.14) and (3.16), and upon equating coefficients of $(ae^{i\theta}, ae^{-i\theta}; q)_k$ on both sides of (3.15) we obtain

$$\frac{(q-1)^2 \lambda}{4} a_k = q^{-k}(1 - q^{k+1})(1 - abq^k)(1 - acq^k)(1 - adq^k) a_{k+1}$$
$$- q^{1-k}(1 - abcdq^{k-1})(1 - q^k) a_k.$$

Thus

$$(3.17) \quad a_{k+1} = \frac{\left[q(1 - abcdq^{k-1})(1 - q^k) + q^k(q-1)^2\lambda/4\right]}{(1 - q^{k+1})(1 - abq^k)(1 - acq^k)(1 - adq^k)} a_k.$$

This shows that $a_{n+1} = 0$ but $a_n \neq 0$ if and only if $\lambda = \lambda_n$. When $\lambda = \lambda_n$ it is straightforward to see that (3.17) implies

$$(3.18) \quad a_k = q^k(q^{-n}, abcdq^{n-1}; q)/(q, ab, ac, ad; q)_q a_0.$$

It can be shown that Theorem 3.4 follows from this.

**4. Remarks.** The Askey-Wilson polynomials form a very important class of orthogonal polynomials which include as special or limiting cases the classical polynomials of Hermite, Laguerre and Jacobi as well as the discrete orthogonal polynomials of Poisson-Charlier, Meixner and Hahn, [5]. They retain most of the structural properties of the traditional polynomials of Jacobi, Hermite and Laguerre. This led Andrews and Askey [3] to redefine classical orthogonal polynomials as the Askey-Wilson polynomials or special or limiting cases of them. The Askey Tableau [14], which indicates how to go from the Askey-Wilson polynomials at the top of the chart down the hierarchy of classical orthogonal polynomials to the Hermite polynomials at the bottom of the chart.

We believe the approach outlined in the first four sections of this work makes it possible to cover the Askey-Wilson polynomials in an introductory course on special functions, for example as a supplement to [17]. The other known proofs of the orthogonality relation of the Askey-Wilson polynomials, [11], [12], [16], except for [6], first evaluate the Askey-Wilson integral, $\int_{-1}^{1} w(x;a,b,c,d)dx$, see (1.6), then use this evaluation and a $_3\phi_2$ summation theorem to establish the orthogonality relation (1.5).

An important feature of this work is that it is essentially self contained and the only results needed are the Jacobi triple product identity (3.12) and the three term recurrence relation for the Askey-Wilson polynomials. In a first course on special functions one can first prove the $q$-binomial theorem

$$(4.1) \qquad \sum_{n=0}^{\infty} \frac{(b;q)_n}{(q;q)_n} z^n = \frac{(bz;q)_\infty}{(z;q)_\infty}, \quad |z| < 1,$$

when $b = q^k, k = 0, 1, \cdots$, by acting repeatedly with the operator $D_q$,

$$(4.2) \qquad (D_q f)(x) := \frac{f(x) - f(qx)}{(1-q)x},$$

on the geometric series

$$\sum_{n=0}^{\infty} z^n = 1/(1-z), \quad |z| < 1.$$

The validity of (4.1) for general $b$ then follows from the identity theorem for analytic functions since both sides of (4.1) are entire functions of $b$. One can then use the same argument to prove the Ramanujan $_1\psi_1$ sum, [7, (5.2.1)]

$$(4.3) \qquad \sum_{-\infty}^{\infty} \frac{(a;q)_n}{(b;q)_n} z^n = \frac{(q, b/a, az, q/az; q)_\infty}{(b, q/a, z, b/az; q)_\infty}, \quad |b/a| < |z| < 1.$$

The Jacobi triple product identity follows if in (4.3) we set $b = 0$, replace $a$ and $z$ by $1/a$ and $az$ respectively, then let $a \to 0$. This proof is due to Ismail, [8] and is reproduced in [1]. Another proof using functional equations is due to Andrews and Askey [2] and is reproduced in [7].

The three term recurrence relation satisfied by the Askey-Wilson polynomials is

$$(4.4) \qquad 2x p_n(x;a,b,c,d|q) = A_n p_{n+1}(x;a,b,c,d|q) + B_n p_n(x;a,b,c,d|q) + C_n p_{n-1}(x;a,b,c,d|q),$$

with

$$(4.5) \quad A_n = \frac{1 - abcdq^{n-1}}{(1 - abcdq^{2n-1})(1 - abcdq^{2n})},$$

$$(4.6) \quad C_n = \frac{(1-q^n)(1-abq^{n-1})(1-acq^{n-1})(1-adq^{n-1})}{(1-abcdq^{2n-2})(1-abcdq^{2n-1})}$$
$$\times (1 - bcq^{n-1})(1 - bdq^{n-1})(1 - cdq^{n-1}),$$

and

$$(4.7) \quad B_n = a + a^{-1} \quad -A_n a^{-1}(1 - abq^n)(1 - acq^n)(1 - adq^n)$$
$$-C_n a[(1 - abq^{n-1})(1 - acq^{n-1})(1 - adq^{n-1})]^{-1}.$$

Wilson [18] has discovered a systematic way of obtaining three term recurrence relations for hypergeometric and basic hypergeometric functions. The recurrence relation (4.4) can also be established by substituting for the $p_n$'s from (1.4), using

$$-2x(ae^{i\theta}, ae^{-i\theta}; q)_k = q^{-k}\left[(ae^{i\theta}, ae^{-i\theta}; q)_{k+1} - (1 + a^2q^{2k})(ae^{i\theta}, ae^{-i\theta}; q)_k\right]$$

then equating coefficients of $(ae^{i\theta}, ae^{-i\theta}; q)_k$ on both sides of (4.4).

As an application we now evaluate the connection coefficients $\{c_{n,j}(a,b)\}$ in (1.18). It is obvious that

$$(4.8) \quad \zeta_j(a, b, aq^{1/2}, bq^{1/2})c_{n,j}(a,b) = \int_{-1}^{1} w(x; aq^{1/2}, bq^{1/2}, aq, bq) \, p_n(x; a, b, aq^{1/2}, bq^{1/2}|q)$$
$$\times p_j(x; aq^{1/2}, bq^{1/2}, aq, bq|q) \, dx.$$

In view of (3.4), (3.5) and (1.12) we have

$$\zeta_j(a, b, aq^{1/2}, bq^{1/2})c_{n,j}(a,b)$$

$$= \left(\frac{q-1}{2}\right) q^{(j-1)/2} \int_{-1}^{1} \left[\mathcal{D}_q w(x; aq, bq, aq^{3/2}, bq^{3/2}) p_{j-1}(x; aq, bq, aq^{3/2}, bq^{3/2}|q)\right]$$
$$\times p_n(x; a, b, aq^{1/2}, bq^{1/2}|q) \, dx$$

$$= q^{(j-n)/2}(1 - q^n)(1 - a^2b^2q^n)\zeta_{j-1}(aq^{1/2}, bq^{1/2}, aq, bq)c_{n-1,j-1}(aq^{1/2}, bq^{1/2}),$$

where we have used (4.8). Therefore (1.7) implies

$$(4.9) \quad c_{n,j}(a,b) = \frac{q^{(j-n)/2}(1-q^n)(1-a^2b^2q^n)}{(1-q^j)(1-a^2b^2q^j)} c_{n-1,j-1}(aq^{1/2}, bq^{1/2}).$$

Repeated application of (4.9) yields

$$(4.10) \quad c_{n,j}(a,b) = \frac{q^{j(j-n)/2}(q;q)_n \left(a^2b^2q^n;q\right)_j}{(q;q)_{n-j}\left(q,a^2b^2q^j;q\right)_j} c_{n-j,0}(aq^{j/2}, bq^{j/2}),$$

which reduces the problem to the evaluation of $c_{m,0}(a,b)$. Note that

$$\frac{w(x; aq^{1/2}, bq^{1/2}, aq, bq)}{w(x; a, b, aq^{1/2}, bq^{1/2})} = (1 - ae^{i\theta})(1 - ae^{-i\theta})(1 - be^{i\theta})(1 - be^{-i\theta})$$
$$= (1 - 2ax + a^2)(1 - 2bx + b^2).$$

Since

$$p_0(x; a, b, aq^{1/2}, bq^{1/2}|q) = 1,$$

and

$$p_1(x; a, b, aq^{1/2}, bq^{1/2}|q) = 2(1 - a^2b^2q)x - (a+b)(1 + q^{1/2})(1 - abq^{1/2}),$$

we find that

$$\frac{w(x; aq^{1/2}, bq^{1/2}, aq, bq)}{w(x; a, b, aq^{1/2}, bq^{1/2})} = \alpha_0 p_0(x; a, b, aq^{1/2}, bq^{1/2}|q) + \alpha_1 p_1(x; a, b, aq^{1/2}, bq^{1/2}|q)$$
$$+ \alpha_2 p_2(x; a, b, aq^{1/2}, bq^{1/2}|q),$$

with

$$(4.11) \quad \alpha_0 = \frac{(bd-1)(bc-1)(ad-1)(ab-1)(ac-1)(abq-1)}{(abcdq-1)(abcd-1)},$$

$$(4.12) \quad \alpha_1 = -\frac{(abq-1)\left(a^2bcdq + ab^2cdq - abdq - aqbc - abd - abc + a + b\right)}{abcdq^2 - 1},$$

and

$$(4.13) \quad \alpha_2 = \frac{ab}{(1 - abcdq)(1 - abcdq^2)}.$$

Thus

$$(4.14) \quad \zeta_0(a, b, aq^{1/2}, bq^{1/2}) c_{m,0}(a,b) = \alpha_0 \zeta_0(a, b, aq^{1/2}, bq^{1/2}) \delta_{m,0}$$

$$+ \alpha_1 \zeta_1(a, b, aq^{1/2}, bq^{1/2}) \delta_{m,1} + \alpha_2 \zeta_2(a, b, aq^{1/2}, bq^{1/2}) \delta_{m,2}.$$

The value of $c_{n,j}(a,b)$ now follows from (4.10)-(4.14) and (1.7).

# 5. The $q$-Sturm-Liouville Operator $T$. Let

$$(5.1) \quad (f, g) := \int_{-1}^{1} f(x)\overline{g(x)}dx.$$

It is evident that (1.12) implies that

$$(5.2) \quad (\mathcal{D}_q f, g) = -(f, \mathcal{D}_q g),$$

holds for all $f$, $g \in H_{1/2}$. This yields the following result.

**Lemma 5.1** *For all $f$, $g \in H_{1/2}$ we have*

$$(5.3) \quad \int_q^{q^{1/2}} \left([\check{f}, \check{g}](x) - [\check{f}, \check{g}](-x)\right) dx = 0,$$

*where*

$$(5.4) \quad [\check{f}, \check{g}](x) = \check{f}(x)\check{g}(q^{-1/2}x) - \check{f}(q^{-1/2}x)\check{g}(x).$$

**Proof**. We first observe that

$$(5.5) \quad \begin{aligned}(\mathcal{D}_q f, g) &= -\frac{iq^{1/2}}{q-1} \int_0^\pi \left\{\check{f}(q^{1/2}e^{i\theta}) - \check{f}(q^{-1/2}e^{i\theta})\right\} \bar{\check{g}}(e^{i\theta}) d\theta \\ &= -\frac{iq^{1/2}}{q-1} \left\{\int_0^\pi \check{f}(q^{1/2}e^{i\theta})\bar{\check{g}}(e^{i\theta})d\theta - \int_{-\pi}^0 \check{f}(q^{1/2}e^{i\theta})\bar{\check{g}}(e^{i\theta})d\theta\right\}.\end{aligned}$$

If $C_2^+$ denotes $\{z = q^{1/2}e^{i\theta} : 0 < \theta < \pi\}$, we have

$$(5.6) \quad \begin{aligned}\int_0^\pi \check{f}(q^{1/2}e^{i\theta})\bar{\check{g}}(e^{i\theta})d\theta &= \int_{C_2^+} \check{f}(q^{1/2}z)\bar{\check{g}}(z)\frac{dz}{iz} \\ &\quad + \left[\int_{-1}^{-q^{1/2}} + \int_{q^{1/2}}^1\right] \check{f}(q^{1/2}x)\bar{\check{g}}(x)\frac{dx}{ix} \\ &= -\int_0^{-\pi} \check{f}(e^{i\theta})\bar{\check{g}}(q^{1/2}e^{i\theta})d\theta + \left[\int_{-1}^{-q^{1/2}} + \int_{q^{1/2}}^1\right] \check{f}(x)\bar{\check{g}}(q^{-1/2}x)\frac{dx}{ix}\end{aligned}$$

since $\bar{\check{g}}(q^{-1/2}e^{-\theta}) = \bar{\check{g}}(q^{1/2}e^{\theta})$. It follows that

$$\begin{aligned}(\mathcal{D}_q f, g) - (f, \mathcal{D}_q g) &= -\frac{iq^{1/2}}{q-1}\{\int_0^\pi [\check{f}(q^{1/2}e^{i\theta})\bar{\check{g}}(e^{i\theta}) + \bar{\check{g}}(q^{1/2}e^{i\theta})\check{f}(e^{i\theta})]d\theta \\ &\quad - \int_{-\pi}^0 [\check{f}(q^{1/2}e^{i\theta})\bar{\check{g}}(e^{i\theta}) + \bar{\check{g}}(q^{1/2}e^{i\theta})\check{f}(e^{i\theta})]d\theta\} \\ &= -\frac{q^{1/2}}{q-1}\left\{\int_{-q^{-1/2}}^{-q} + \int_q^{q^{-1/2}}\right\}[\check{f}\check{g}](x)\frac{dx}{x} \\ &= \frac{q^{1/2}}{1-q}\int_q^{q^{1/2}}([\check{f}\check{g}](x) - [\check{f}\check{g}](-x))\frac{dx}{x}\end{aligned}$$

whence (5.3).

Let $\mathcal{H}_w$ denote the weighted space $L^2(-1, 1; w(x)dx)$ with inner product

$$(f, g)_w := \int_{-1}^{1} f(x)\overline{g(x)}w(x)dx, \quad \|f\|_w := (f, f)_w^{1/2}$$

and let $T$ be defined by

$$Tf(x) := Mf(x)$$

for $f$ in $H_1$, where

(5.7) $\quad (Mf)(x) = -\dfrac{1}{w(x)}\mathcal{D}_q\left(p\mathcal{D}_q f\right)(x).$

We shall assume that $p$ and $w$ are positive on $(-1, 1)$ and also satisfy

(5.8) $\quad$ (i) $\quad p(x)/\sqrt{1-x^2} \in H_{1/2}$, $1/p \in L(-1, 1)$,

$\quad\quad\quad$ (ii) $\quad w(x) \in L(-1, 1), 1/w \in L(-1, 1; \dfrac{dx}{(1-x^2)})$.

The expression $Mf$ is therefore defined for $f \in H_1$, and the operator $T$ acts in $\mathcal{H}_w$. Furthermore, the domain $H_1$ of $T$ is dense in $\mathcal{H}_w$ since it contains all polynomials.

**Lemma 5.2** *The operator $T$ is symmetric in $\mathcal{H}_w$ and $T \geq 0$.*

**Proof.** We infer from (5.2) that for all $f, g \in H_1$,

$$\begin{aligned}(5.9) \quad (Tf, f)_w &= -(\mathcal{D}_q[p\mathcal{D}_q f], f) = (p\mathcal{D}_q f, \mathcal{D}_q f) \\ &= \int_{-1}^{1} p(x)|\mathcal{D}_q f(x)|^2 dx,\end{aligned}$$

whence the Lemma.

Let $\mathcal{Q}(T)$ denote the form domain of $T$ and $\tilde{T}$ its Friedrichs extension. Recall that $\mathcal{Q}(T)$ is the completion of $H_1$ with respect to $\|.\|_\mathcal{Q}$, where

(5.10) $\quad \|f\|_\mathcal{Q}^2 := \int_{-1}^{1} p(x)|\mathcal{D}_q f|^2 dx + \|f\|_w^2,$

and if $(.,.)_\mathcal{Q}$ denotes the inner product on $\mathcal{Q}(T)$, then for all $f \in \mathcal{D}(\tilde{T})$ and $g \in \mathcal{Q}(T)$,

(5.11) $\quad (f, g)_\mathcal{Q} = ([\tilde{T} + I]f, g)_w.$

where $I$ is the identity on $\mathcal{H}_w$. We have that $f \in \mathcal{Q}(T)$ if and only if there exists a sequence $\{f_n\} \subset H_1$ such that $\|f - f_n\|_\mathcal{Q} \to 0$; hence $\|f - f_n\|_w \to 0$ and $\{\mathcal{D}_q f_n\}$ is a Cauchy sequence in $L^2(-1, 1; p(x)dx)$, with limit $F$ say. From (5.8) and (5.2) it follows that for $\phi \in H_{1/2}$,

$$(5.12) \quad \int_{-1}^{1} F(x)\phi(x)dx = \lim_{n \to \infty} \int_{-1}^{1} (\mathcal{D}_q f_n)(x)\phi(x)dx = -\lim_{n \to \infty} \int_{-1}^{1} f_n(x)\mathcal{D}_q \phi(x)dx,$$
$$= -\int_{-1}^{1} f(x)\mathcal{D}_q \phi(x)dx + \mathcal{O}\left(\|f - f_n\|_w [\int_{-1}^{1} |\mathcal{D}_q \phi(x)|^2 \frac{dx}{w(x)}]^{1/2}\right)$$
$$= -\int_{-1}^{1} f(x)\mathcal{D}_q \phi(x)dx.$$

Thus, in analogy with distributional derivatives, we shall say that $F = \mathcal{D}_q f$ in the generalized sense. Note that this proves that $F$ is unique up to a function that vanishes almost everywhere, so different Cauchy sequences $f_n$ give the same $\mathcal{D}_q f$.

We conclude that the norm on $\mathcal{Q}(T)$ is defined by (5.10) with $\mathcal{D}_q f$ now understood in the generalized sense. Also, it follows in a standard way that

$$(5.13) \quad \mathcal{D}(T^*) = \{f : f, Mf \in \mathcal{H}_w\}, \quad T^*f = Mf,$$

$$(5.14) \quad \mathcal{D}(\tilde{T}) = \mathcal{Q}(T) \cap \mathcal{D}(T^*)$$
$$= \{f : p^{1/2}\mathcal{D}_q f \in L^2(-1, 1), \ Mf; \in \mathcal{H}_w\}.$$

If $T$ is the operator in the Askey-Wilson case, that is

$$w(x) = w(x; a, b, c, d), \quad p(x) = w(x; aq^{1/2}, bq^{1/2}, cq^{1/2}, dq^{1/2}),$$

then the Askey-Wilson polynomials satisfy (1.15), that is

$$Tp_n(x; a, b, c, d|q) = -\lambda_n p_n(x; a, b, c, d|q).$$

Since $p_n(x; a, b, c, d|q)$ is of degree $n$, $(n = 0, 1, ...)$ and the polynomials are dense in $\mathcal{H}_w$, it follows that $\{p_n(x; a, b, c, d|q)\}$ forms a basis for $\mathcal{H}_w$. Hence $T$ has a selfadjoint closure $\overline{T} = \tilde{T}$ and $\tilde{T}$ has a discrete spectrum consisting of the eigenvalues $\lambda_n$ in (1.16), $n = 0, 1, ....$ If $f \in \mathcal{Q}(T)$ then setting $p_n(x) \equiv p_n(x; a, b, c, d|q)$, we have

$$(f, p_n)_\mathcal{Q} = (f, [T + I]p_n)_w = (\lambda_n + 1)(f, p_n)_w$$

and, in particular, from (1.5) with $\zeta_n \equiv \zeta_n(a, b, c, d)$,

$$(p_m, p_n)_\mathcal{Q} = (\lambda_n + 1)\zeta_n \delta_{m,n}.$$

It follows that $e_n = p_n/\sqrt{\zeta_n(\lambda_n + 1)}$, $(n = 0, 1, ...)$ is an orthonomal basis for $\mathcal{Q}(T)$. Thus $f \in \mathcal{Q}(T)$ if and only if

$$(5.15) \quad \sum_{n=0}^{\infty} |(f, e_n)_{\mathcal{Q}}|^2 = \sum_{n=0}^{\infty} \frac{\lambda_n + 1}{\zeta_n} |(f, p_n)_w|^2 < \infty,$$

$f \in \mathcal{D}(\tilde{T})$ if and only if

$$(5.16) \quad \sum_{n=0}^{\infty} |([\tilde{T} + I]f, \frac{p_n}{\zeta_n^{1/2}})_w|^2 = \sum_{n=0}^{\infty} \frac{(\lambda_n + 1)^2}{\zeta_n} |(f, p_n)_w|^2 < \infty$$

and, for $f \in \mathcal{D}(\tilde{T})$

$$(5.17) \quad \tilde{T}f = \sum_{n=0}^{\infty} (\tilde{T}f, \frac{p_n}{\zeta_n^{1/2}})_w \frac{p_n}{\zeta_n^{1/2}} = \sum_{n=0}^{\infty} \frac{\lambda_n}{\zeta_n} (f, p_n)_w p_n.$$

A special case of interest is that of the Chebyshev polynomials of the first kind where

$$(5.18) \quad a = -b = 1, c = -d = q^{1/2}, \quad w(x) = (1 - x^2)^{-1/2}, \quad p(x) = 4(1 - x^2)^{1/2},$$

$p_n(x) = T_n(x)$, and $T_n(\cos\theta) = \cos n\theta$. Here we have, in the notation (1.9), with $x = \cos\theta$ and $z = e^{i\theta}$, that

$$\sqrt{1 - x^2} \mathcal{D}_q f(x) = \frac{1}{i(q^{1/2} - q^{-1/2})} F(x),$$

where

$$F(x) = \breve{F}(z) = \breve{f}(q^{1/2}z) - \breve{f}(q^{-1/2}z).$$

Thus

$$(Mf)(x) := -\sqrt{1 - x^2} \mathcal{D}_q \{4\sqrt{1 - x^2} \mathcal{D}_q f(x)\} = -4\{\frac{\breve{F}(q^{1/2}e^{i\theta}) - \breve{F}(q^{-1/2}e^{i\theta})}{(i[q^{1/2} - q^{-1/2}])^2}\}$$

$$= \frac{4}{(q^{1/2} - q^{-1/2})^2} \{\breve{f}(qe^{i\theta}) - 2\breve{f}(e^{i\theta}) + \breve{f}(q^{-1}e^{i\theta}).$$

For $f(x) = T_n(x)$

$$\breve{f}(z) = (z^n + z^{-n})/2$$

from which it follows that

$$\sqrt{1-x^2}\mathcal{D}_q\{\sqrt{1-x^2}\mathcal{D}_q T_n(x)\} = -\frac{\{q^{n/2}-q^{-n/2}\}^2}{\{q^{1/2}-q^{-1/2}\}^2}T_n(x).$$

**Acknowledgments**. We thank Richard Askey for reminding us of reference [13] and for providing a further simplification of the evaluation of the connection coefficients in (1.18). We are grateful to Dennis Stanton for sending us several corrections. Last but not least we thank the referee for a detailed report offering many helpful criticisms and suggestions to improve the presentation of this work. This work was done while M. Ismail was visiting University of Wales College at Cardiff and he gratefully acknowledges the hospitailty of the School of Mathematics and the Department of Computing Mathematics.

Department of Computing Mathematics, University of Wales College at Cardiff, Cardiff, Wales CF2 4YN
School of Mathematics, University of Wales College at Cardiff, Cardiff, Wales CF4 2YN
Department of Mathematics, University of South Florida, Tampa, FL 33620, U. S. A.